\documentclass[12pt]{amsart}
\usepackage[cp1251]{inputenc}
\usepackage[english]{babel}
\usepackage{amsmath}
\usepackage{amssymb}
\usepackage{amsfonts}

\begin{document}


\title[ON INFINITE MATRICES]
      {ON INFINITE MATRICES}

  \author{A.V. Vasilyev}
\address{Chair of Applied Mathematics and Computer Modeling\\
    Belgorod State National Research University\\
         Pobedy street 85, Belgorod 308015, Russia}

        \email{alexvassel@gmail.com}

\author{V.B. Vasilyev}
\address{Chair of Applied Mathematics and Computer Modeling\\
    Belgorod State National Research University\\
         Pobedy street 85, Belgorod 308015, Russia}

        \email{vbv57@inbox.ru}

 \author{A.B. Kamanda Bongay}
\address{Chair of  Applied Mathematics and Computer Modeling\\
 Belgorod State National Research University\\
 Studencheskaya street 14, Belgorod 308007, Russia}

\email{159720@bsu.edu.ru}

\begin{abstract}
We consider linear bounded operators acting in Banach spaces with a basis, such operators can be represented by an infinite matrix. We prove that for an invertible operator there exists a sequence of invertible finite-dimensional operators so that the family of norms of their inverses is uniformly bounded. It leads to the fact that solutions of finite-dimensional equations converge to the solution of initial operator equation with infinite-dimensional matrix.
\end{abstract}

\maketitle

\section{Introduction}

Discrete equations is a very important mathematical object. This is related to computer calculations which help us to find numerical solution if we don't know its analytical expression.
Discrete equations can appear via difference schemes \cite{S} or difference potentials \cite{R}, or discrete convolutions \cite{GF}. The latter is more interesting for us because we try to develop discrete theory for pseudo-differential equations based on ideas and methods \cite{E,V0}. Certain realization of these ideas and methods is presented in authors' papers \cite{VV6,V1,VV1,VV2,VT,VVT,VVM}. All mentioned papers are related to a solvability problem for discrete pseudo-differential equations. Such equations are roughly speaking infinite systems of linear algebraic equations, and for numerical solution we need to approximate these infinite systems by certain finite systems. In such cases, they used the reduction method.

This reduction method was developed in \cite{GF} for abstract situation and for different classes of operators. Some results were obtained in papers \cite{K,KS} for general operators and discrete convolutions. But these papers don't give an answer to the question if arbitrary invertible operator admits the reduction method, assuming the operator is presented by infinite matrix. In this paper, we will prove this assertion. .

\section{Infinite matrices}

Let $X$ be a Banach space with standard basis ${\bf e}_i=(0,\dots,0,\underbrace{1}_{i},0,\dots),$ and we consider an infinite system of linear algebraic equations with the matrix $A=(a_{ij})_{i,j=1}^{\infty}$; this matrix is a representation of linear bounded operator $A$ in the space $X$:
\[
A: X\rightarrow X.
\]

Let's introduce the following equation in the space $X$
\begin{equation}\label{1}
A{\bf x}={\bf y},~~~{\bf y}\in X.
\end{equation}

Let $P_n$ be a projector on a linear span of vectors  ${\bf e}_i, i=1,\dots,n$; this linear span will be denoted by $X_n$. Then we put $A_n=P_nAP_n$ so that $A_n: X_n\rightarrow X_n$ and we write the truncated equation
\begin{equation}\label{2}
A_n{\bf x}_n=P_n{\bf y}
\end{equation}
in the vector space $X_n$. Thus, the operator $A_n$ is represented by the matrix $(a_{ij})_{i,j=1}^n$. Obviously, the sequence of operators $A_n$ strongly converges to $A$, i.e. $\forall{\bf x}\in X,\lim\limits_{n\to\infty}A_n{\bf x}=A{\bf x}$.

We will give here one auxiliary result which will help us to obtain main theorem.

{\bf Lemma.} {\it If a certain ${\bf x}\in X,\lim\limits_{n\to\infty}A_n{\bf x}=A{\bf x}$, and there is the sequence $\{{\bf x}_n\}_{n=1}^{\infty}\subset X$ such that $\lim\limits_{n\to\infty}{\bf x}_n={\bf x}$ then $\lim\limits_{n\to\infty}A_n{\bf x}_n=A{\bf x}$.
}

\begin{proof}
Indeed, we have
\[
||A{\bf x}-A_n{\bf x}_n||\leq||A{\bf x}-A_n{\bf x}||+||A_n{\bf x}-A_n{\bf x}_n||\leq
\]
\[
\leq||A{\bf x}-A_n{\bf x}||+||A||\cdot||{\bf x}-{\bf x}_n||.
\]

Both summands $||A{\bf x}-A_n{\bf x}||$ and $||A||\cdot||{\bf x}-{\bf x}_n||$ tend to zero according to assumptions of Lemma, and the proof is completed.
\end{proof}

\section{Main result}

The following assertion is called usually the {\it "reduction method"}.

\textbf{Theorem.} {\it If the inverse bounded operator $A^{-1}: X\rightarrow X$ exists then the following assertions are valid:

1) starting from a certain $N, \forall n\geq N,$ the operators $A_n: X_n\rightarrow X_n$ are invertible;

2) we have the estimate
\[
||A_n^{-1}||\leq C,
\]
with constant $C$ non-depending on $n$;

 3) the solution ${\bf x}_n$ to the equation $(2)$ converges to the solution ${\bf x}$ of the equation $(1)$ under $n\to\infty$.
}

\begin{proof}
We use the proof by contradiction applying the theory of finite systems of linear algebraic equations. Namely, the Cramer's rule asserts that the operator $A_n$ with the matrix $(a_{ij})_{i,j=1}^n$ will be invertible in the space $X_n$ iff $\det A_n\neq 0$.

{\it First step.} Let's note that there are two possibilities for considered situation: either
starting from a certain $N, \forall n\geq N,$ the operators $A_n: X_n\rightarrow X_n$ are invertible or there is a subsequence $A_{n_k}$ of non-invertible operators. If the first situation is valid then we have the needed assertion. That's why we assume that the second situation is realized. So, we have a sequence of non-invertible operators $A_{n_k},\,k\to\infty$. We will show that this assumption leads to a contradiction.

If operators $A_{n_k}$ are invertible then $\det A_{n_k}=0$. But then there exists such a matrix $(a_{ij})_{i,j=1}^{n_k}$ with a certain non-zero minor so that all minor of bigger order are vanishing, otherwise all matrices $A_{n_k}$ will be null-matrices. We will assume (without loss of generality) that this minor is related to the matrix
$$
A_{m_k=}
\begin{pmatrix}
a_{11}&a_{12}&\cdots a_{1m_k}\\
a_{21}&a_{22}&\cdots a_{2m_k}\\
\cdots&\cdots&\cdots\\
a_{m_k1}&a_{m_k2}&\cdots a_{m_km_k}\\
\end{pmatrix},\quad
~~~~~~~~~~~~~~~~~~~~~m_k<n_k,
$$
and $\det A_{m_k}\neq 0$. This implies that the homogeneous system
\[
A_{n_k}{\bf x}_{n_k}=0
\]
has non-trivial solutions. Now we will describe their structure.

Let's denote by $Y_{m_k}$ the space $X_{n_k}\ominus X_{m_k}$ so that
\[
X_{n_k}=X_{m_k}\oplus Y_{m_k},
\]
and the representation
\[
{\bf x}_{n_k}={\bf x}_{m_k}+{\bf y}_{m_k},~~~{\bf x}_{m_k}\in X_{m_k},~~~{\bf y}_{m_k}\in Y_{m_k},
\]
is unique for arbitrary ${\bf x}_{n_k}\in X_{n_k}$.

We introduce the rectangular  $(m_k\times n_k)$-matrix $B_{m_k}$ of the following type
\[
B_{m_k=}
\begin{pmatrix}
a_{1m_k+1}&a_{1m_k+2}&\cdots a_{1n_k}\\
a_{2m_k+1}&a_{2m_k+2}&\cdots a_{2n_k}\\
\cdots&\cdots&\cdots\\
a_{m_km_k+1}&a_{m_km_k+2}&\cdots a_{m_kn_k}\\
\end{pmatrix}.
\]

Such operator $B_{m_k}$ is uniformly bounded as an operator $Y_{n_k-m_k}\rightarrow X_{m_k}$ because
\[
||B_{m_k}y_{m_k}||\leq||Ay_{m_k}||.
\]
Therefore, we have the following property
\[
{\bf x}_{m_k}=-A^{-1}_{m_k}B_{m_k}{\bf y}_{m_k}.
\]

If we will transfer to a limit then we will see that the operator $A$ has infinite-dimensional kernel, thus, it is non-invertible, and we have contradiction.

{\it Second step.}
According to the first step we have that starting from a certain $N, \forall n\geq N,$ the operators $A_n: X_n\rightarrow X_n$ are invertible. Let's assume that the sequence $||A_n^{-1}||$ is unbounded. It means that there are sequences $\{{\bf x}_{n}\}_{n=1}^{\infty}\subset X_n$ and $\{c_{n}\}_{n=1}^{\infty},\,c_n>0$ such that
\[
||A^{-1}_n{\bf x}_n||\geq c_n||{\bf x}_n||,~~~~~\lim\limits_{n\to\infty}c_n=\infty.
\]
If we put ${\bf x}'_n={\bf x}_n/||{\bf x}_n||$ so that $||{\bf x}'_n||=1$ then we can write
\begin{equation}\label{3}
||A^{-1}_n{\bf x}'_n||\geq c_n.
\end{equation}

Let ${\bf y}_n=A^{-1}_n{\bf x}'_n$. According to \eqref{3} we have ${\bf y}_n\to\infty$. Then
\[
A_n{\bf y}_n={\bf x}'_n,~~~~~||A_n{\bf y}_n||=1.
\]
We put ${\bf y}'_n={\bf y}_n/||{\bf y}_n||,\,||{\bf y}'_n||=1$ and then
\[
||A_n{\bf y}'_n||=1/||{\bf y}_n||\rightarrow 0,\,n\to\infty.
\]

Thus, we have the sequence $\{{\bf y}'_n\}_{n=1}^{\infty},\,||{\bf y}'_n||=1$ such that
\begin{equation}\label{4}
||A_n{\bf y}'_n||\rightarrow 0,\,n\to\infty.
\end{equation}

Let's consider $A_n{\bf y}'_n-A{\bf y}'_n$. This is a vector of the following type
\[
\begin{pmatrix}
0\\
\cdots\\
0\\
\sum\limits_{k=1}^{n}a_{n+1k}y'_k\\
\sum\limits_{k=1}^{n}a_{n+2k}y'_k\\
\cdots
\end{pmatrix},
\]
where first $n$ coordinates are zero. It seems That this vector tends to zero under $n\to\infty$. But there is a counterexample, the basis $\{{\bf e}_k\}_{k=1}^{\infty}$. Nevertheless, we will find a contradiction using another way.

Let's consider the sequence of operators $A_n$ more carefully. Obviously, according to \eqref{4} we have
\[
\inf\limits_{||{\bf y}||=1}||A_n{\bf y}||=\alpha_n,~~~\lim\limits_{n\to\infty}\alpha_n=0.
\]

Further, for an arbitrary ${\bf y}\in X,\,||{\bf y}||=1$ we have
\[
\left|||A_n{\bf y}||-||A{\bf y}||\right|\leq||A_n{\bf y}-A{\bf y}||,
\]
and then
\[
||A{\bf y}||\leq||A_n{\bf y}||+||A_n{\bf y}-A{\bf y}||.
\]
Thus,
\[
\inf\limits_{||{\bf y}||=1}||A{\bf y}||\leq\inf\limits_{||{\bf y}||=1}||A_n{\bf y}||+||A_n{\bf y}-A{\bf y}||\leq\alpha_n+||A_n{\bf y}-A{\bf y}||
\]
for all ${\bf y}\in X$. Fix ${\bf y}\in X$. Given $\varepsilon>0$ we can find such $N\in\mathbb N$ that $\forall n\geq N$ we have
\[
\alpha_n
\]
so that
\[
\inf\limits_{||{\bf y}||=1}||A{\bf y}||<\varepsilon,
\]
and we conclude
\begin{equation}\label{5}
\inf\limits_{||{\bf y}||=1}||A{\bf y}||=0.
\end{equation}

The equality \eqref{5} implies that there is the sequence $\{{\bf z}_k\}_{k=1}^{\infty},\,||{\bf z}||_k=1$ such that
\[
\lim\limits_{k\to\infty}A{\bf z}_k=0.
\]
But the operator $A$ is invertible, and then $\lim\limits_{k\to\infty}{\bf z}_k=0$. The latter assertion is a contradiction.

{\it Third step.} This step is related to a convergence. We have
\[
{\bf x}=A^{-1}{\bf y}
\]
and
\[
{\bf x}_n=A_n^{-1}P_n{\bf y}
\]
We will denote $P_n{\bf y}={\bf y}_n$ and estimate ${\bf x}-{\bf x}_n$. Let's write
\[
{\bf x}-{\bf x}_n=A^{-1}{\bf y}-A_n^{-1}{\bf y}_n=
\]
and thus,
\[
||{\bf x}-{\bf x}_n||=||A^{-1}{\bf y}-A_n^{-1}{\bf y}_n||.
\]
Now we can apply Lemma with $A^{-1}$ and $A^{-1}_n$ instead of
$A$ and $A_n$, ant Theorem is proved.

\end{proof}

{\bf Remark.}
May be such a result exists in  mathematical literature but the authors have no appropriate information.

\section*{Conclusion}

This studying is very important for our studying discrete pseudo-differential equations and related discrete boundary value problems. As a rule such problems lead to infinite systems of linear algebraic equations, and we need a verification for change the infinite system by finite one. Moreover, pseudo-differential operators are defined in Fourier images, and now is not clear what approach is more effective from computational point of view, original space or its Fourier image.


\end{document}